\documentclass{amsart}
\usepackage[T1]{fontenc}

\title{Solution of the Basel problem using the Feynman integral trick}
\author{F. L. Freitas}
\email{felipelopesfreitas@gmail.com}
\date{\today}

\begin{document}

\begin{abstract}
Euler's solution in 1734 of the Basel problem, which asks for a closed form expression for the sum of the reciprocals of all perfect squares, is one of the most celebrated results of mathematical analysis. In the modern era, numerous proofs of it have been produced, each emphasizing a different style of calculation, as a way of testing the power of some demonstration method. It's often thought that solutions using calculus need to involve advanced contour integration techniques or geometric coordinate transformations. We show that this is not the case, as the result can be derived by analyzing basic properties of a particular one-dimensional integral, and as such, can be obtained with techniques typical of regular calculus tests and math competitions.
\end{abstract}

\maketitle

\section{Introduction}

Solving the Basel problem amounts to proving the following identity:

\begin{equation}
\sum_{n=1}^\infty \frac{1}{n^2} = 1 + \frac{1}{2^2} + \frac{1}{3^2} + \frac{1}{4^2} + \dots = \frac{\pi^2}{6}
\end{equation}

The problem was first solved by Euler in 1734\cite{ayoub1974}. Ever since his famous original proof, which relied on a generalized Vieta formula for the roots of the function $\frac{\sin x}{x}$, numerous other proofs have appeared\cite{chapman2003}, most emphasizing a different method to attack the problem.

Of all of these, the most elementary is perhaps that of Apostol\cite{apostol1983}, where the problem is reduced to evaluating a two-dimensional integral:

\begin{equation}
\sum_{n=1}^\infty \frac{1}{n^2} = \int_0^1\int_0^1 \frac{\mathrm{d}x\mathrm{d}y}{1-xy}
\end{equation}

Apostol noted that, if we perform the 2D coordinate transformation $(u,v)=((x+y)/2,(x-y)/2)$, the integral can be evaluated exactly, and the $\pi^2/6$ answer emerges.

We give a proof similar in spirit to the one of Apostol, but arguably simpler, because it uses only a 1D integral and instead of an \emph{ad hoc} coordinate transformation, it 
employs differentiation under the integral sign, a technique popularized by physicist Richard Feynman\cite{feynman2010}.

\section{Result}

We attempt to evaluate the integral:

\begin{equation}
\int_0^\infty \ln\left(1 + e^{-2x}\right)\mathrm{d}x
\end{equation}

However, instead of evaluating it as is, we solve a more general problem, where the integral has an extra parameter $\alpha$:

\begin{equation}\label{eq:integ}
I(\alpha) = \int_0^\infty \ln\left(1 + \alpha e^{-x} + e^{-2x}\right)\mathrm{d}x
\end{equation}

Although it looks like an extra complication, the idea is that the dependency of the integral on $\alpha$ can be very simple, and can be deduced from a differential equation in the parameter $\alpha$, which might be easier to solve than the original integral. If we differentiate, we obtain:

\begin{equation}
\frac{\mathrm{d}I}{\mathrm{d}\alpha} = \int_0^\infty \frac{e^{-x}}{1 + \alpha e^{-x} + e^{-2x}}\mathrm{d}x
\end{equation}

This integral can be simplified by changing variables to $u=e^{-x}$, where we get $\mathrm{d}u=-e^{-x}\mathrm{d}x$. We can get rid of the minus sign by swapping the integration limits. Renaming $u \to x$, we can write:

\begin{equation}
\frac{\mathrm{d}I}{\mathrm{d}\alpha} = \int_0^1 \frac{1}{1 + \alpha x + x^2}\mathrm{d}x
\end{equation}

Since are interested in the range $-2<\alpha<2$, the integrand can be evaluated with standard techniques, by completing the square in the denominator. The answer is:

\begin{equation}
\frac{\mathrm{d}I}{\mathrm{d}\alpha} = \left.\left[ \frac{2}{\sqrt{4-\alpha^2}}\arctan\left(\frac{\alpha+2x}{\sqrt{4-\alpha^2}}\right)\right] \right|_{x=0}^1
\end{equation}

Substituting the limits:

\begin{equation}
\frac{\mathrm{d}I}{\mathrm{d}\alpha} = \frac{2}{\sqrt{4-\alpha^2}}\left(\arctan\left(\frac{\alpha+2}{\sqrt{4-\alpha^2}}\right)-\arctan\left(\frac{\alpha}{\sqrt{4-\alpha^2}}\right)\right)
\end{equation}

The expression can be simplified with the $\arctan$ identity:

\begin{equation}
\arctan(x) - \arctan(y) = \arctan\left(\frac{x-y}{1+xy}\right)
\end{equation}

We obtain:

\begin{equation}
\frac{\mathrm{d}I}{\mathrm{d}\alpha} = \frac{2}{\sqrt{4-\alpha^2}}\arctan\left(\sqrt{\frac{2-\alpha}{2+\alpha}}\right)
\end{equation}

Although the expression looks formidable, it can be cast into an elementary form by substituting $\alpha=2\cos u$. The argument of $\arctan$ simplifies with the following trigonometric identity:

\begin{equation}
\sqrt{\frac{1-\cos u}{1+\cos u}} = \sqrt{\frac{2\sin^2(u/2)}{2\cos^2(u/2)}} = \tan\left(\frac{u}{2}\right)
\end{equation}

And the end result is:

\begin{equation}\label{eq:intalpha}
I(\alpha) = -\frac{1}{2}\arccos\left(\frac{\alpha}{2}\right)^2 + c
\end{equation}

We only need to determine the integration constant $c$. To do that, notice that we have a factorization in the integrand when $\alpha=2$:

\begin{equation}
I(2) = \int_0^\infty \ln\left(1 + 2e^{-x} + e^{-2x}\right)\mathrm{d}x = 2\int_0^\infty \ln\left(1 + e^{-x}\right)\mathrm{d}x
\end{equation}

By applying the transformation $x\to2x$, we can reduce it to the integral at $\alpha=0$:

\begin{equation}
2\int_0^\infty \ln\left(1 + e^{-x}\right)\mathrm{d}x = 4\int_0^\infty \ln\left(1 + e^{-2x}\right)\mathrm{d}x = 4I(0)
\end{equation}

The relation $I(2)=4I(0)$ allows us to determine $c$. We have:

\begin{equation}
I(2) = c \qquad I(0) = -\frac{\pi^2}{8} + c
\end{equation}

Therefore:

\begin{equation}
I(2) = 4I(0) \Rightarrow c = -\frac{\pi^2}{2} + 4c \Rightarrow c = \frac{\pi^2}{6}
\end{equation}

Now we are able to determine \eqref{eq:integ} for any value of $\alpha$. In particular, we can calculate it for $\alpha=-2$. The corresponding integral is:

\begin{equation}
I(-2) = \int_0^\infty \ln\left(1 - 2e^{-x} + e^{-2x}\right)\mathrm{d}x = 2\int_0^\infty \ln\left(1 - e^{-x}\right)\mathrm{d}x
\end{equation}

Substituting $\alpha=-2$ into \eqref{eq:intalpha}, we obtain:

\begin{equation}\label{eq:baselint}
I(-2) = -\frac{\pi^2}{3} \Rightarrow \int_0^\infty \ln\left(1-e^{-x}\right)\mathrm{d}x = -\frac{\pi^2}{6}
\end{equation}

From the geometric series expansion, we have:

\begin{equation}\label{eq:geom}
\frac{1}{1-x} = 1 + x + x^2 + x^3 + \dots = \sum_{n=0}^\infty x^n
\end{equation}

Integrating \eqref{eq:geom} term by term we find:

\begin{equation}
-\ln(1-x) = x + \frac{x^2}{2} + \frac{x^3}{3} + \dots = \sum_{n=1}^\infty \frac{x^n}{n}
\end{equation}

Substituting $e^{-x}$ and integrating we find:

\begin{equation}\label{eq:baselexp}
-\int_0^\infty \ln(1-e^{-x})\mathrm{d}x = \sum_{n=1}^\infty \int_0^\infty \frac{e^{-nx}}{n}\mathrm{d}x = \sum_{n=1}^\infty \frac{1}{n^2}
\end{equation}

Together, equations \eqref{eq:baselint} and \eqref{eq:baselexp} imply our desired result:

\begin{equation}
\sum_{n=1}^\infty \frac{1}{n^2}= \frac{\pi^2}{6}
\end{equation}

\section{Conclusion}

It has been shown that the solution to the classic Basel problem can be found in a straightforward manner using the Feynman integral trick. The method applied here is quite unusual, because in most applications of differentiation under the integral sign, one obtains the integration constant by going to a limit where the integral is trivial.

Here, the determination of the constant is done by looking at two configurations of the parametrized integral related by a symmetry. This is a quite interesting technique, which might have applications in other integrals. The difficulty of applying the Feynman trick rests, as always, on how to introduce the parameter to simplify the differentiated integral. However, like many of these tricks, when it's seen once, it's much easier to apply it in new settings.

In summary, we have provided a very simple proof of the Basel problem which relies on the Feynman trick in a one-dimensional integral. Although hard to find, after the parameter is introduced, the proof follows naturally, which means such a technique might be useful in situations where a solution has to be produced in tight time constraints, as is often the case in math competitions.

\bibliographystyle{ieeetr}
\bibliography{basel}

\end{document}